
\documentstyle{amsppt}
\baselineskip18pt
\magnification=\magstep1
\pagewidth{30pc}
\pageheight{45pc}
\hyphenation{co-deter-min-ant co-deter-min-ants pa-ra-met-rised
pre-print pro-pa-gat-ing pro-pa-gate
fel-low-ship Cox-et-er dis-trib-ut-ive}
\def\leaderfill{\leaders\hbox to 1em{\hss.\hss}\hfill}
\def\A{{\Cal A}}
\def\C{{\Cal C}}

\def\idest{i.e.,\ }

\def\a{{a}}
\def\be{{b}}
\def\g{{c}}

\def\d{{\delta}}

\def\e{{\varepsilon}}

\def\l{{\lambda}}

\def\bc{{\bold c}}

\def\b0{\text{\bf 0}}

\def\ra{{\ \longrightarrow \ }}
\def\sra{{\rightarrow}}

\def\ds{\displaystyle}

\def\coker{\text{\rm \, coker }}

\def\lan{{\langle}}
\def\ran{{\rangle}}

\def\zed{{\Bbb Z}}

\def\Im{\text{\rm Im}}

\def\pd{\partial}
\def\boxit#1{\vbox{\hrule\hbox{\vrule \kern3pt
\vbox{\kern3pt\hbox{#1}\kern3pt}\kern3pt\vrule}\hrule}}
\def\rabbit{\vbox{\hbox{\kern0pt
\vbox{\kern0pt{\hbox{---}}\kern3.5pt}}}}

\def\tableau#1{
        \hbox {
                \hskip -10pt plus0pt minus0pt
                \raise\baselineskip\hbox{
                \offinterlineskip
                \hbox{#1}}
                \hskip0.25em
        }
}

\def\tabCol#1{
\hbox{\vtop{\hrule
\halign{\strut\vrule\hskip0.5em##\hskip0.5em\hfill\vrule\cr\lower0pt
\hbox\bgroup$#1$\egroup \cr}
\hrule
} } \hskip -10.5pt plus0pt minus0pt}

\def\CR{
        $\egroup\cr
        \noalign{\hrule}
        \lower0pt\hbox\bgroup$
}

\def\mapright#1{
    \mathop{\longrightarrow}\limits^{#1}}
\def\mapdown#1{\Big\downarrow
     \rlap{$\vcenter{\hbox{$\scriptstyle#1$}}$}}

\def\CD#1{
   $$
   \def\normalbaselines{\baselineskip20pt \lineskip3pt \lineskiplimit3pt }
    \matrix #1 \endmatrix
   $$
}

\def\blank#1#2{
\hbox to #1{\hfill \vbox to #2{\vfill}}
}


\def\strut{\vrule height10pt depth5pt width0pt}

\topmatter
\title Acyclic heaps of pieces, I
\endtitle

\author R.M. Green \endauthor
\affil 
Department of Mathematics and Statistics\\ Lancaster University\\
Lancaster LA1 4YF\\ England\\
{\it  E-mail:} r.m.green\@lancaster.ac.uk\\
\endaffil

\abstract
Heaps of pieces were introduced by Viennot and have applications to algebraic 
combinatorics, theoretical computer science and statistical physics.  In this
paper, we how certain combinatorial properties of heaps studied by Fan 
and by Stembridge are closely related to the properties of a certain
linear map $\pd_E$ associated to a heap $E$.  We examine the
relationship between $\pd_E$ and $\pd_F$ when $F$ is a subheap of $E$.
This approach allows neat statements and proofs 
of results on certain associative algebras (generalized Temperley--Lieb 
algebras) that are otherwise tricky to prove.  The key to the proof is to 
interpret the structure constants of the aforementioned algebras 
in terms of the maps $\pd$.
\endabstract

\thanks
\noindent 2000 {\it Mathematics Subject Classification.} 06A11.
\endthanks


\endtopmatter

\centerline{\bf To appear in the Journal of Algebraic Combinatorics}


\def\hker#1{\ker \pd_{#1}}
\def\hcoker#1{\coker \pd_{#1}}

\head Introduction \endhead

A heap is an isomorphism class of labelled posets 
satisfying certain axioms.  Heaps have a wide
variety of applications, notably to parallelism in computer
science, but also to statistical physics and algebraic combinatorics.  Many 
of these applications are discussed by Viennot in \cite{{\bf 17}}.

One of the oldest results in algebraic topology states that if $k$ is a field
and $G$ is a finite, connected, directed graph with vertex set $V(G)$ and edge
set $E(G)$, then the linear map $\pd : kE(G) \ra kV(G)$ sending the edge 
$(v_i \ra v_j)$ to $v_i - v_j$ has image of codimension $1$ in $k V(G)$.

In this paper, we will introduce and study an analogue of the above
situation for heaps.
The definitions are completely general but we are
particularly interested in heaps arising from fully commutative elements
in Coxeter groups as studied by Stembridge in \cite{{\bf 15}}.
We will show how this theory may be applied to obtain neater 
proofs of results on certain associative algebras.

We now summarise the main results of the paper for easy reference.  Section
1 sets up the basic definitions of heaps, including our definition of 
the map $\pd$ in Definition 1.2.1.  A heap $E$ will be called acyclic if
$\pd_E$ is injective, and $E$ will be called strongly acyclic
if it and all its maximal subheaps are acyclic. 
The main body of theory is developed
in Section 2.  We look at two combinatorial properties of heaps, properties
P1 and P2, and show how they are related to the linear notions of
being acyclic and strongly acyclic.  In favourable circumstances, 
property P1 is equivalent to being acyclic (Theorem 2.4.4) and 
property P2 is equivalent to being strongly acyclic (Theorem 2.4.2).  

We are interested in how $\ker \pd_E$ is related to $\ker \pd_F$ when
$F$ is a subheap of $E$.  The Deletion Lemma (Theorem 2.1.1)
shows that $\ker \pd$ changes in dimension by at most $1$ when 
an element is removed from the heap.  In \S2.3, we look at certain specific
kinds of subheaps for which $\ker \pd$ is exactly predictable from those
of the original heap.  These constructions exhibit a close relationship
between the map $\pd$ and a certain quotient of the heap monoid algebra whose
structure constants can be described in terms of $\pd$
(Theorem 3.2.3).  We apply the theory here to give simpler proofs of
certain results involving generalized Temperley--Lieb algebras; one such
result is given in Proposition 3.4.2, and more are given in \S4.1.  In \S4.2
we sketch the relationship between $\ker \pd$ and certain diagram
calculi, and we conclude with some questions in \S4.3.

Although we were led to this theory by questions about generalized 
Temperley--Lieb algebras, it is hoped that the theory in this paper will be
of independent interest.


\head 1. Heaps \endhead

In \S1.1, we introduce the basic properties of heaps.  We will tend to 
follow Viennot's notation \cite{{\bf 17}}.  We give our definition of
$\pd$ in \S1.2.

\subhead 1.1 Basic definitions \endsubhead

\definition{Definition 1.1.1}
Let $P$ be a set equipped with a symmetric and reflexive binary relation
$\C$.  The elements of $P$ are called {\it (basic) pieces}, and the relation
$\C$ is called the {\it concurrency relation}.

A {\it labelled heap} with pieces in $P$ is a triple $(E, \leq, \e)$ 
where $(E, \leq)$ is a finite (possibly empty)
partially ordered set with order relation denoted
by $\leq$ and $\e$ is a map $\e : E \ra P$ satisfying the following two
axioms. 

\item{1.}{For every $\a, \be \in E$ such that $\e(\a) \ \C \ \e(\be)$, 
$\a$ and $\be$ are comparable in the order $\leq$.}

\item{2.}{The order relation $\leq$ is the transitive closure of the
relation $\leq_\C$ such that for all $\a, \be \in E$, $\a \leq_\C \be$ 
if and only if both $\a \leq \be$ and $\e(\a) \ \C \ \e(\be)$.}
\enddefinition

We will sometimes express the relation $\a \leq \be$ by saying that 
``$\a$ is {\it above} $\be$''.  The terms {\it minimal} and 
{\it maximal} applied to the elements of the labelled 
heap refer to minimality (respectively, maximality) with respect to $\leq$.

Parts 1 and 2 of Definition 1.1.1  correspond to axioms (i) and (ii'')
respectively in Viennot's paper.

\example{Example 1.1.2}
Let $P = \{1, 2, 3\}$ and, for $x, y \in P$, define $x \ \C \ y$ 
if and only if 
$|x - y| \leq 1$.  Let $E = \{a, b, c, d, e\}$ partially ordered by
extension of the (covering) relations $a \leq c$, $b \leq c$, $c \leq d$,
$c \leq e$.  Define the map $\e$ by the conditions $\e(a) = \e(d) = 1$,
$\e(c) = 2$ and $\e(b) = \e(e) = 3$.  Then $(E, \leq, \e)$ can easily be
checked to satisfy the axioms of Definition 1.1.1 and it is a labelled heap.
The minimal elements are $a$ and $b$, and the maximal elements are $d$ and $e$.
\endexample

\definition{Definition 1.1.3}
Let $(E, \leq, \e)$ and $(E', \leq', \e')$ be two labelled 
heaps with pieces in $P$ and with the same concurrency relation, $\C$.  
An isomorphism $\phi : E \ra E'$ of posets is said to be an 
{\it isomorphism of labelled posets} if $\e = \e' \circ \phi$.

A {\it heap} of pieces in $P$ with concurrency relation $\C$ is a labelled
heap (Definition 1.1.1) defined up to labelled poset isomorphism.
The set of such heaps is denoted by $H(P, \C)$.  We denote the heap 
corresponding to the labelled heap $(E, \leq, \e)$ by $[E, \leq, \e]$.
\enddefinition

We will sometimes abuse language and speak of the underlying set of a heap,
when what is meant is the underlying set of one of its representatives.

\definition{Definition 1.1.4}
Let $(E, \leq, \e)$ be a labelled heap with pieces in $P$ and let $F$ be a 
subset of $E$.  
Let $\e'$ be the restriction of $\e$ to $F$.  Let ${\Cal R}$ be the relation
defined on $F$ by $\a \ {\Cal R} \ \be$ if and only if $\a \leq \be$ and
$\e(\a) \ \C \ \e(\be)$.  Let $\leq'$ be the transitive closure of ${\Cal R}$.
Then $(F, \leq', \e')$ is a labelled heap with pieces in $P$.  The heap 
$[F, \leq', \e']$ is called a {\it subheap} of $[E, \leq, \e]$.
\enddefinition

We will often implicitly use the fact that a subheap is determined by its
set of vertices and the heap it comes from.

\definition{Definition 1.1.5}
The {\it concurrency graph} associated to the class of heaps $H(P, \C)$ is
the graph whose vertices are the elements of $P$ and for which there is an
edge from $v \in P$ to $w \in P$ if and only if $v \ne w$ and $v \ \C \ w$.
\enddefinition

\definition{Definition 1.1.6}
Let $E = [E, \leq_E, \e]$ and $F = [F, \leq_F, \e']$ be two heaps in 
$H(P, \C)$.
We define the heap $G = [G, \leq_G, \e''] = E \circ F$ of $H(P, \C)$ 
(which we call the superposition of $E$ over $F$) as follows.

\item{1.}{The underlying set $G$ is the disjoint union of $E$ and $F$.}
\item{2.}{The labelling map $\e''$ is the unique map $\e'' : G \ra P$ whose
restriction to $E$ (respectively, $F$) is $\e$ (respectively, $\e'$).}
\item{3.}{The order relation $\leq_G$ is the transitive closure of the
relation ${\Cal R}$ on $G$, where $\a \ {\Cal R} \ \be$ if
and only if one of the following three conditions holds:
\item{(i)}{$\a, \be \in E$ and $\a \leq_E \be$;}
\item{(ii)}{$\a, \be \in F$ and $\a \leq_F \be$;}
\item{(iii)}{$\a \in E, \ \be \in F$ and $\e(\a) \ \C \ \e'(\be)$.}
}
\enddefinition

\remark{Remark 1.1.7}
Definition 1.1.6 can easily be shown to be sound (see \cite{{\bf 17}, \S2}).
It is immediate from the construction that $E$ and $F$ are subheaps of 
$E \circ F$.  Note that Viennot calls $E \circ F$ ``the superposition 
of $F$ over $E$''.  

As in \cite{{\bf 17}}, we will write $\a \circ E$ and $E \circ \a$ for
$\{\a\} \circ E$ and $E \circ \{ \a \}$, respectively.  Note that $\a \circ E$
and $\be \circ E$ are equal as heaps if $\e(\a) = \e(\be)$.
\endremark

\definition{Definition 1.1.8}
A {\it trivial heap} is a heap $[E, \leq, \e]$ for which the order relation
$\leq$ is trivial, meaning that no element of $E$ is above any other.
\enddefinition

\subhead 1.2 The map $\pd$ \endsubhead

We can now introduce our analogue of the graph theoretic phenomenon mentioned
in the introduction; this is the central definition of this paper.  
Throughout \S1.2, we let $[E, \leq, \e]$ be a heap
in the set $H(P, \C)$ with pieces in $P$ and concurrency relation $\C$.  
We also fix a field, $k$.

\definition{Definition 1.2.1}
Let $V_0$ be the set of elements of $[E, \leq, \e]$, \idest the set of 
elements of (a representative of) the underlying poset, $E$.
We call the elements of $V_0$ {\it vertices} and  denote their $k$-span 
by $C_0$.

Let $V_1$ be the set of all pairs $(x, y) \in E \times E$ with $x < y$ and 
$\e(x) = \e(y)$ such that there is no element $z$ for which we have
both $\e(x) = \e(z) = \e(y)$ and $x < z < y$.  
We call the elements of $V_1$ {\it edges} and denote their $k$-span by $C_1$.

For all other integers $i \in \zed \backslash \{0, 1\}$, we define $C_i = 0$.

The $k$-linear map $\pd  = \pd_E : C_1 \ra C_0$ is defined by its 
effect on the edges as follows: $$
\pd : (x, y) \mapsto \sum_{{x < w < y} \atop {\e(w) \ \C \ \e(x)}} w
.$$
\enddefinition

\remark{Remark 1.2.2}
Note that, in the sum of Definition 1.2.1, we have $\e(x) = \e(y)$ but 
it is not possible for $\e(w) = \e(x)$ 
because of the conditions imposed on the edge $(x, y)$.
\endremark

\example{Example 1.2.3}
Consider the heap $[E, \leq, \e]$ arising from the labelled heap of Example
1.1.2.  In this case, $C_1 = \{(a, d), (b, e)\}$ and 
$C_0 = \{a, b, c, d, e\}$.  We have $\pd((a, d)) = c$ and $\pd((b, e)) = c$.
It follows that $\hker{E}$ is of dimension $1$ spanned by 
$[(a, d) - (b, e)]$ and that $\hcoker{E}$ is of dimension $4$, spanned
by the images of $a, b, d, e$.
\endexample

The following simple lemma explains the relationship between the
dimensions of $\ker \pd$ and $\coker \pd$.

\proclaim{Lemma 1.2.4}
Let $[E, \leq, \e]$ be a heap in $H(P, \C)$.  Then $$
\dim \hcoker{E} - \dim \hker{E}
$$ is $|\e(E)|$, the number of
different labels occurring on vertices of the heap.
\endproclaim

\demo{Proof}
It is clear from Definition 1.2.1 that $$
\dim \hcoker{E} - \dim \hker{E} = \dim C_0 - \dim C_1
.$$  For each $p \in P$, let $n_p(E) = |\e^{-1}(p)|$; note that $\e^{-1}(p)$ 
is a chain in $E$.  The definition of $V_1$
shows that the number of edges
$(x, y) \in V_1$ with $\e(x) = \e(y) = p$ is $\max(0, \ n_p(E) - 1)$.  This
shows that $|V_1| = |E| - |\e(E)|$.  
The lemma follows from the fact that $|V_0| = |E|$.
\qed\enddemo

It will be convenient to use the following definitions in the sequel.

\definition{Definition 1.2.5}
Let $E = [E, \leq_E, \e]$ be a heap in $H(P, \C)$ and let $k$ be a field.  
If $v \in E$, we let $E(v) = [E(v), \leq_{E(v)}, \e']$ be the subheap of $E$ 
obtained by defining $E(v) = E \backslash \{v\}$.

We say $E$ is {\it acyclic} if $\hker{E} = 0$.  We say $E$ is {\it strongly
acyclic} if $E$ is acyclic and $E(v)$ is acyclic for all $v \in E$.  We say
$v$ is an {\it image vertex} of $E$ if $v \in \Im(\pd_E)$.
\enddefinition

A large portion of this paper will be concerned with the following two 
problems and their applications.

\proclaim{Problem 1.2.6}
Can we characterize and classify the acyclic (respectively, strongly 
acyclic) heaps in $H(P, \C)$?
\endproclaim

\proclaim{Problem 1.2.7}
If $E$ is a heap in $H(P, \C)$ and $F$ is a subheap of $E$, how is $\hker{F}$
related to $\hker{E}$?
\endproclaim

\head 2. Properties of the map $\pd$ \endhead

In \S2 we prove some results about the map $\pd$ that hold in a general 
context.  In \S2.1, we will prove a deletion lemma that relates 
$\hker{E}$ to $\hker{E(v)}$, in the notation of Definition 1.2.5.  
In \S2.2, we exhibit some general relationships between combinatorial and
linear properties of heaps.
In \S2.3, we show how to ``contract'' a heap to a
simpler one for which $\ker \pd$ is very similar; this turns out to be
a very useful proof technique.
We introduce the notion of a regular class
of heaps $H(P, \C)$ in \S2.4; the heaps of these classes have particularly
tractable properties.  

Our motivation behind these results is to apply them to certain associative
algebras in \S3, but the results of \S2 are related to each other in
intriguing ways that shed light on problems 1.2.6 and 1.2.7.

As before, we will fix a set $H(P, \C)$ and a field, $k$.

\subhead 2.1 The deletion lemma \endsubhead

The next theorem is the main result of \S2.1.  It is very useful for certain 
applications, including our proof of Proposition 3.4.2.

\proclaim{Theorem 2.1.1 (Deletion Lemma)}
Let $[E, \leq, \e]$ be a nonempty heap in \newline 
$H(P, \C)$, and let $v \in E$.  Then $
|\dim \hker{E} - \dim \hker{E(v)}| \leq 1
.$
\endproclaim

To prove this, it is convenient to define the following vector spaces and
maps.

\definition{Definition 2.1.2}  Fix $v$ as in Theorem 2.1.1.
\item{(i)}{Let $A_1$ be the space $C_1$ associated with the heap $E$.  
Let $A_0$ be the quotient of
the space $C_0$ associated to $E$ by the $1$-dimensional subspace
$\lan v \ran$.
Let $\pd_A : A_1 \ra A_0$ be the composition of the map
$\pd$ associated to $E$ with the natural epimorphism.}
\item{(ii)}
{Let $B_1$ be the space $C_1$ associated to $E(v)$, let $B_0$ be the space
$C_0$ associated to $E(v)$ and let $\pd_B$ be the map $\pd$ associated to
$E(v)$.}
\item{(iii)}
{We define a $k$-linear 
map $f_1 : B_1 \ra A_1$ by its effect on the edges of $B$
as follows.  If $(x, y) \in B_1$ is such that $\e(x) = \e(y) = \e(v)$
and $x < v < y$ in $A$, we define $f_1((x, y)) = (x, v) + (v, y)$.  Otherwise,
we define $f_1((x, y)) = (x, y)$.  (This is well defined by the definition
of subheaps and Definition 1.2.1.)}
\item{(iv)}
{Let $f_0$ be the obvious $k$-isomorphism from $B_0$ to $A_0$.}
\enddefinition

\proclaim{Lemma 2.1.3}
In the notation of Definition 2.1.2, the following diagram commutes and has
exact rows.

\CD{
0 & \mapright{} & B_1 & \mapright{f_1} & A_1 & \mapright{} & \coker f_1 & 
\mapright{} & 0 \cr
& & \mapdown{\pd_B} & & \mapdown{\pd_A} & & \mapdown{} & \cr
0 & \mapright{} & B_0 & \mapright{f_0} & A_0 & \mapright{} & 0 & 
\mapright{} & 0 \cr
}
\endproclaim

\demo{Proof}
Injectivity of $f_1$ follows from the fact that the
images of basis elements in $B_1$ under $f_1$ have disjoint supports in $A_1$.
(Recall that the fibres of $\e$ are totally ordered.)

The rest of the claim 
follows from Definition 1.2.1, Definition 2.1.2 and a diagram chase.
\qed\enddemo

The Snake Lemma immediately gives the following result.

\proclaim{Lemma 2.1.4}
Maintain the above notation.  There is an exact sequence $$
0 \ra \hker{E(v)} \ra \ker \pd_A \mapright{} \coker f_1 \mapright{}
\hcoker{E(v)} \ra \coker \pd_A \ra 0.
$$ \qed
\endproclaim

In order to extract more information from Lemma 2.1.4, we need to know
whether $v$ is an image vertex or not (see Definition 1.2.5).
\vfill\eject
\proclaim{Lemma 2.1.5}
Maintain the above notation.
\item{\rm (i)}{If $v$ is an image vertex then
$$\dim \ker \pd_A = \dim \hker{E} + 1$$ and  
$$\dim \coker \pd_A = \dim \hcoker{E}.$$}
\item{\rm (ii)}{If $v$ is not an image vertex then 
$$\dim \ker \pd_A = \dim \hker{E}$$ and  
$$\dim \coker \pd_A = \dim \hcoker{E} - 1.$$}
\endproclaim

\demo{Proof}
This is an exercise in linear algebra, using the definition of $\pd_A$.
\qed\enddemo

\proclaim{Lemma 2.1.6}
If $v$ is the unique vertex $x$ of $E$ with $\e(v) = \e(x)$, then 
$f_1$ is surjective.  Otherwise, $\dim \coker f_1 = 1$.
\endproclaim

\demo{Proof}
If $v$ satisfies the uniqueness property above then there can be no 
edges $(x, y)$ in $A_1$ (or $B_1$) with either $x = v$ or $y = v$.  It
follows that $f_1$ is an isomorphism in this case.

On the other hand, if $v$ does not satisfy the uniqueness property then
the argument of Lemma 1.2.4 shows that $\dim A_1 = \dim B_1 + 1$, and the
lemma follows.
\qed\enddemo

\demo{Proof of Theorem 2.1.1}
Suppose that $v \in \Im(\pd_E)$; this is case (i) of Lemma 2.1.5.
Since $\dim \coker f_1 \leq 1$ by Lemma 2.1.6, Lemma 2.1.4 shows that
$\dim \hker{E(v)} - \dim \hker{E}$ is equal to $0$ or $1$.

Suppose now that $v \not\in \Im(\pd_E)$.  A similar argument based on
case (ii) of Lemma 2.1.5 shows that
$\dim \hker{E(v)} - \dim \hker{E}$ is equal to $0$ or $-1$.
\qed\enddemo

\subhead 2.2 Heaps with additional properties \endsubhead

The two properties P1 and P2 introduced in
this section have combinatorial definitions, but as we shall see, they are
related to properties of the map $\pd$ and they cast some light on Problem
1.2.6.

\definition{Definition 2.2.1 (Property P1)}
Let $E = [E, \leq, \e] \in H(P, \C)$ be a heap.
We write $E(\a) \prec^+ E$ (respectively, $E(\a) \prec^- E$) 
if $\a$ is a maximal (respectively, minimal) vertex of $E$ and there
exists a maximal (respectively, minimal) vertex $\be$ of $E(\a)$ with
$\e(\be) \ne \e(\a)$ such that $\be$ is not 
maximal (respectively, minimal) in $E$.  We write $E(\a) \prec E$ if either
$E(\a) \prec^+ E$ or $E(\a) \prec^- E$.

If there is a (possibly trivial) sequence $
E_1 \prec E_2 \prec \cdots \prec E
$ of heaps in $H(P, \C)$ where $E_1$ is a trivial heap, we say that 
the heap $E$ is {\it dismantlable} or that $E$ has {property P1}.
\enddefinition

\example{Example 2.2.2}
The heap $E$ arising from Example 1.1.2 does not have property P1, but its 
subheap $E(a)$ does: consider the sequence of subheaps $$
\{ d, e \} \prec \{ c, d, e \} \prec \{ b, c, d, e \}
.$$
\endexample

In \S3.4, we shall exploit the relationship between property P1 and
Fan's notion of
left and right cancellability \cite{{\bf 5}, Definition 4.2.4}.  We avoid
the term ``cancellability'' in this paper because of possible confusion 
with the use of this term in the theory of monoids.

\proclaim{Proposition 2.2.3}
A dismantlable heap is acyclic.
\endproclaim

\demo{Proof}
Let $E$ be a dismantlable heap and let $$
E_1 \prec E_2 \prec \cdots \prec E_l = E
$$ be a chain of heaps with $E_1$ trivial.  The proof is by induction on $l$.
If $l = 1$, the claim is clear because a trivial heap has no edges.

For the inductive step, we treat the case where $E_{l-1} \prec^- E_l$; the 
other case is similar.  Let $\a$ be a minimal element of $E_l$ and let
$\be$ be a minimal vertex of $E_{l-1}$ that is not minimal in
$E_l$ with $\e(\be) \ne \e(\a)$.  Suppose $\hker{E_l} \ne 0$ and let $$
\sum \l_i e_i
$$ be a nontrivial element of $\ker \pd$.  Since $\hker{E_{l-1}} = 0$, one
of the edges $e_i$ must involve the vertex $\a$; let us write $e_i = (\a, \g)$
as $\a$ is minimal.  (Note that $\g \ne \be$ as $\e(\g) = \e(\a) \ne \e(\be)$.)
Since $\be$ is minimal in $E_{l-1}$, the vertex $\be$
cannot occur with nonzero coefficient in $\pd(e_j)$ for any edge
$e_j$ with $j \ne i$.  However, $\be$ occurs with coefficient $1$ in 
$\pd(e_i)$, contrary to hypothesis.
\qed\enddemo

\remark{Remark 2.2.4}
The converse of Proposition 2.2.3 is false in general.  
Consider the class of heaps
$H(P, \C)$ for which the concurrency graph is a square whose corners (the 
elements of $P$) are consecutively labelled $p_1, p_2, p_3, p_4$.  Let $E
\in H(P, \C)$ be a labelled heap whose underlying set is 
$\{a_1, a_2, a_3, a_4\}$ and $\e(a_i) = p_i$ for all $i$.  
Let $a_i < a_j$ whenever 
$i$ is odd and $j$ is even.  The heap $[E, \leq, \e]$ is acyclic but not
dismantlable.
\endremark

The next property is modelled on Stembridge's characterization of full
commutativity in \cite{{\bf 15}, Proposition 2.3}.  The term ``convex
chain'' for a heap has its natural meaning: a chain $$
x_1 < x_2 < \cdots < x_t
$$ of vertices 
in a heap is convex if and only if whenever $x_i < y < x_j$ for some
$y$, the vertex $y$ is an element of the chain.

\definition{Definition 2.2.5 (Property P2)}
We say a heap $E = [E, \leq, \e] \in H(P, \C)$
has property P2 if it contains no convex chains of the form $x < y < z$ or
$x < z$ with $\e(x) = \e(z)$ in either case.  
\enddefinition

\example{Example 2.2.6}
The heap arising from Example 1.1.2 does not have property P2.  Although there
are no chains of the form $x < z$ with $\e(x) = \e(z)$, the chains $a < c < d$
and $b < c < e$ each violate the other requirement.
\endexample

\proclaim{Proposition 2.2.7}
A strongly acyclic heap has property P2.
\endproclaim

\demo{Proof}
Let $E = [E, \leq, \e]$ be a heap that fails property P2.  If $E$ contains
an convex chain of the form $x < z$ with $\e(x) = \e(z)$ then $(x, z)$ is
an edge in $E$ and $\pd((x, z)) = 0$, meaning that $E$ is not acyclic.
The other possibility is that $E$ contains a convex chain $x < y < z$
with $\e(x) = \e(z)$.  In this case, the subheap $E(y)$ contains an edge
$(x, z)$ with $\pd((x, z)) = 0$, and $E(y)$ is not acyclic, meaning that
$E$ is not strongly acyclic.
\qed\enddemo

\remark{Remark 2.2.8}
The converse of Proposition 2.2.7 is false in general.  Consider the set
of heaps $H(P, \C)$ defined in Remark 2.2.4.  Let $E$ be a labelled heap
with elements $\{a_1, a_2, a_3, a_4, a_5, a_6\}$ labelled by the elements
$p_1, p_3, p_2, p_4, p_1, p_3$ respectively.
The Hasse diagram of $(E, \leq)$ is shown in Figure 1.
This gives a heap with edges $e_1 = (a_1, a_5)$ and $e_2 = (a_2, a_6)$.  
We have $\pd(e_1) = a_3 + a_4$ and $\pd(e_2) = a_3 + a_4$.  This means that
the heap arising from $E$ is not acyclic, although the reader may easily
check that the heap has property P2.
\endremark

\topcaption{Figure 1} Hasse diagram for the labelled heap $(E, \leq)$ of
Remark 2.2.8
\endcaption
\centerline{
\hbox to 0.791in{
\vbox to 0.875in{\vfill
        \includegraphics{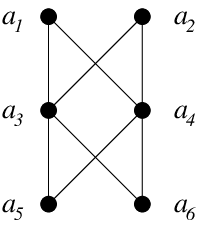}
}
\hfill}
}

The following lemma is similar to \cite{{\bf 15}, Lemma 3.1}.  (See Remark 1.1.7
for the notation.)

\proclaim{Lemma 2.2.9}
Let $E = [E, \leq, \e]$ be a heap in $H(P, \C)$.  If $E$ has property
P2 but $\a \circ E$ does not, then either there is a minimal element $\g$
of $E$ with $\e(\g) = \e(\a)$ or there is a convex chain $\a < \be < d$ in
$\a \circ E$ with $\e(\a) = \e(d) \ne \e(\be)$.
\endproclaim

\demo{Proof}
This is an easy consequence of the definition of property P2.
\qed\enddemo
\vfill\eject
\subhead 2.3 More on convex chains \endsubhead

In \S2.3, we give a precise answer to Problem 1.2.7 in certain special cases
by showing that the subheaps $F$ of a given heap $E$ arising from a certain 
construction are such that $\ker \pd_F$ is predictable from $\ker \pd_E$.
To define this construction, we need the concept
of contraction along a convex chain.

\definition{Definition 2.3.1}
Let $E = [E, \leq, \e]$ be a heap of $H(P, \C)$, and let $$
\bc = (x_1, x_2, \ldots, x_i) : x_1 < x_2 < \cdots < x_i
$$ be a chain in $E$.  We say $\bc$ is {\it balanced} if $\e(x_1) = \e(x_i)$.
If $\bc$ is a balanced convex chain, we define the heap $E/\bc$ to be the 
subheap of $E$ obtained by omitting the vertices $x_2, x_3, \ldots x_i$.
We call the heap $E/\bc$ the {\it contraction of $E$ along $\bc$}, and the
number $i$ is called the {\it length} of the chain.
\enddefinition

\remark{Remark 2.3.2}
If, in the above definition, we had omitted vertices 
$x_1, x_2, \ldots, x_{i-1}$ instead, we would have obtained the same heap:
the fact that $\bc$ is balanced makes the two corresponding labelled heaps
isomorphic.
\endremark

\example{Example 2.3.3}
Take the heap arising from Example 1.1.2.  The chain $b < c < e$ is balanced
and convex, and contraction along this chain gives the subheap with vertices
$\{a, b, d\}$.
\endexample

\proclaim{Lemma 2.3.4}
Let $E = [E, \leq, \e]$ be a heap and let $\bc$ be a balanced convex chain of
$E$ of length $2$.  Then $\dim(\hker{E}) = \dim(\hker{E/\bc}) + 1$.
\endproclaim

\demo{Proof}
Let us write $B = E/\bc$ for notational convenience.  The chain $\bc$ consists
of (a covering pair of) elements $x < y$ with $\e(x) = \e(y)$.  

We denote by $B_1$ (respectively, $E_1$) the span of the edges in $B$
(respectively, $E$), and we denote by $B_0$
(respectively, $E_0$) the span of the vertices in $B$ (respectively, $E$).
Let $g_1$ be the map from $B_1$ to $E_1$ defined by its effect on edges as
follows: $$
g_1((a, b)) = \cases
(y, b) & \text{ if } a = x,\cr
(a, b) & \text{ otherwise.}
\endcases
$$  Let $g_0$ be the map from $B_0$ to $E_0$ defined by its effect on
vertices as follows: $$
g_0(v) = \cases
x + y & \text{ if } v = x,\cr
v & \text{ otherwise.}
\endcases
$$  Definition 1.2.1 implies that the left square in the diagram
\CD{
0 & \mapright{} & B_1 & \mapright{g_1} & E_1 & \mapright{} & \coker g_1 & 
\mapright{} & 0 \cr
& & \mapdown{\pd_B} & & \mapdown{\pd_E} & & \mapdown{0} & \cr
0 & \mapright{} & B_0 & \mapright{g_0} & E_0 & \mapright{} & \coker g_0 & 
\mapright{} & 0 \cr
} 

\noindent commutes.  The right square commutes 
because $\dim \coker g_1 = 1$ and the edge $(x, y) \in E_1$ lies in
$\ker \pd_E \backslash \Im(g_1)$.  The same argument shows that the connecting
homomorphism $z$ in the exact sequence $$
0 \, \sra \, \hker{B} \, \sra \, \hker{E} \mapright{} \coker g_1 \mapright{z}
\hcoker{B} \mapright{} \hcoker{E} \, \sra \, \coker g_0 \, \sra \, 0
$$ provided by the Snake Lemma is zero, and the result follows.
\qed\enddemo

\proclaim{Lemma 2.3.5}
Let $E = [E, \leq, \e]$ be a heap and let $\bc$ be a balanced convex chain 
$x < y < z$ such that $\e(x) \ne \e(y)$.  
Then $\dim \hker{E} = \dim \hker{E/\bc}$.
\endproclaim

\demo{Proof}
Note that $y$ is necessarily an image vertex as it is the image of the
edge $(x, z)$ under $\pd$.  We will write $B = E/\bc$ for notational
convenience.

Suppose that $y$ is the unique vertex $c$ in $E$ with $\e(y) = \e(c)$.
Lemmas 2.1.4, 2.1.5 (i) and 2.1.6 then imply that $$
\dim \hker{E(y)} = \dim \hker{E} + 1
,$$ and Lemma 2.3.4 shows that $$
\dim \hker{E(y)} = \dim \hker{E/\bc} + 1
,$$ because $E/\bc = E(y)/\bc'$ where $\bc'$ is the balanced convex 
chain $x < z$ of $E(y)$.
This proves the claim in this case, and we may now suppose that there exists
a vertex $c \ne y$ with $\e(c) = \e(y)$.  It is either possible to choose 
$c$ so that $(c, y)$ is an edge or to choose $c$ so that $(y, c)$ is an edge.
We will assume the former case is possible; the other case follows mutatis
mutandis.

We denote by $B_1$ (respectively, $B_0$) the span of the edges (respectively,
vertices) of $B$.  Let $A_1$ be the span of all edges of $E$ except $(x, z)$,
and let $A_0$ be the quotient of the space $C_0$ (as in Definition 1.2.1) 
by the $1$-dimensional subspace $\lan y \ran$.  The differential $\pd_E$
induces a $k$-linear map $\pd_A$ from $A_1$ to $A_0$.  Let $h_1$ be the map 
from $B_1$ to $A_1$ defined by its effect on edges as follows: $$
h_1((a, b)) = \cases
(a, y) + (y, b) & \text{ if } (a, y) \text{ and } (y, b) \text{ are edges in }
E,\cr
(z, b) & \text{ if } a = x,\cr
(a, b) & \text{ otherwise.}
\endcases
$$  Let $h_0$ be the map from $B_0$ to $A_0$ defined by its effect on vertices
as follows: $$
h_0(v) = \cases
x + z + \lan y \ran & \text{ if } v = x,\cr
v + \lan y \ran & \text{ otherwise.}
\endcases
$$  A routine check using the definitions shows that the left square in the
diagram 
\CD{
0 & \mapright{} & B_1 & \mapright{h_1} & A_1 & \mapright{} & \coker h_1 & 
\mapright{} & 0 \cr
& & \mapdown{\pd_B} & & \mapdown{\pd_A} & & \mapdown{s} & \cr
0 & \mapright{} & B_0 & \mapright{h_0} & A_0 & \mapright{} & \coker h_0 & 
\mapright{} & 0 \cr
} 

\noindent commutes. 
This induces a map $s$ making the diagram commute.  The map
$s$ is nonzero, since $\coker h_1$ is spanned by the image of $(c, y)$,
$x$ occurs in the support of $\pd((c, y))$ and $z$ does not.
Since $\coker h_0$ is spanned by the image of $x$, $s$ is an isomorphism and
the Snake Lemma applied to this situation gives the exact sequence $$
0 \ra \hker{B} \ra \ker \pd_A \ra 0 \ra \hcoker{B} \ra \coker \pd_A \ra 0
.$$  A linear algebra argument shows that $\dim \ker \pd_A = \dim \ker \pd_E$,
because $$\pd_E((x, z)) = y.$$  The conclusion follows.
\qed\enddemo

\remark{Remark 2.3.6}
It follows easily from the definition of property P2 that if $E$ is a heap
that does not have property P2 then it is possible to contract a 
balanced convex chain in $E$ of one of the types given in lemmas 2.3.4 
and 2.3.5.  This is the main point of the above two results.
\endremark

\subhead 2.4 Regular classes of heaps \endsubhead

We now introduce the notion of a regular class of heaps, and show that in 
a regular class of heaps, the converses to propositions 2.2.3 and 2.2.7 hold.
We shall look at some examples of regular classes of heaps in \S3.4.

\definition{Definition 2.4.1}
A class of heaps $H(P, \C)$ is said to be {\it regular} if any heap of
$H(P, \C)$ with property P2 also has property P1.
\enddefinition

The counterexamples in remarks 2.2.4 and 2.2.8 come from classes of heaps
that are not regular.

\proclaim{Theorem 2.4.2}
Suppose that $H(P, \C)$ is a regular class of heaps.  Let 
$E = [E, \leq, \e']$ be a heap of $H(P, \C)$ with property P2.  Then
\item{\rm (i)}{$E$ is strongly acyclic;}
\item{\rm (ii)}
{if $\hker{\a \circ E} \ne 0$ (respectively, $\hker{E \circ \a} \ne 0$)
then there is a minimal (respectively, maximal) vertex $\g$ of $E$ such that
$\e(\g) = \e(\a) \in P$, where $\e$ is the map associated to the heap 
$\a \circ E$ (respectively, $E \circ \a$).}

In particular, property P2 and the property of being strongly acyclic
coincide for heaps of $H(P, \C)$.
\endproclaim

\demo{Proof}
The last claim is immediate from (i) and Proposition 2.2.7.

The proof of (i) and (ii) is by induction on $n$, the number of vertices 
in the heap $E$.
If $n = 0$, there is nothing to prove.  If $E$ is nonempty but trivial, 
claims (i) and (ii) follow easily, and this deals with the
case $n = 1$.  Let $P(l)$ be the statement ``claim (i) holds when $E$ is 
a heap with $l$ vertices and properties P1 and P2'', and let $Q(l)$ be the
statement ``claim (ii) holds when $E$ is a heap with $l-1$ vertices and 
properties P1 and P2''.  We will be done if we can show that $P(l) \Rightarrow
Q(l+1)$ and $(P(l) \wedge Q(l)) \Rightarrow P(l+1)$.

Suppose $P(l)$ holds, that $E$ is a heap with $l$ vertices and that
$\hker{\a \circ E} \ne 0$.  (We omit consideration of the case 
$\hker{E \circ \a} \ne 0$, which is similar.)  It cannot be the case that
$\a \circ E$ has property P2, because if it did, it would have property P1
by assumption and would be acyclic by Proposition 2.2.3.  By Lemma 2.2.9,
we see that if the statement $Q(l+1)$ fails for the heap $\a \circ E$,
there must be a convex chain $\bc = \a < \be < d$ in $\a \circ E$ with $\e(\a)
= \e(d) \ne \e(\be)$.  Remark 2.3.2 shows that $(\a \circ E)/\bc = E(\be)$,
and by Lemma 2.3.5 we have $$
\dim \hker{\a \circ E} = \dim \hker{E(\be)}
.$$  Since $E$ is strongly acyclic, the right hand side is zero and we have a
contradiction, proving $Q(l+1)$.

Now suppose $P(l)$ and $Q(l)$ hold and that $E$ has $l + 1$ vertices.  
We may assume that $E$ is not trivial.  Because
$E$ has property P1, we have $E' \prec E$ for some heap $E$.  We deal with
the case where $E' \prec^- E$, the other case being similar.  In this case
we have $E = \a \circ E'$ for some vertex $\a$.  Since $E$ has property P1
by assumption, it is acyclic by Proposition 2.2.3 and 
there is a minimal element $\be \in E'$ (with $\e(\be) \ne \e(\a)$)
that is not minimal in $E$.  

Suppose also that $E$ is not strongly acyclic.  Then 
then there is an element $v \in E$ such that
$\hker{E(v)} \ne 0$.  We cannot have $v = \a$ because $E' = E(\a)$ 
inherits property P2 from $E$ and is therefore acyclic.  Suppose $v \ne \be$
and let $\ds{\left[\sum \l_i e_i \right]}$ 
be a nontrivial element of $\hker{E(v)}$, where
the $e_i$ are edges in $E(v)$.  It is not possible for any of
the edges $e_i$ to involve the vertex $\a$, because $\be \ne v$ would 
occur with coefficient $1$ in the image of any edge $(\a, \g)$ but would
not occur in the image of any other edge as $\be$ is minimal in $E'$.
This means that $\ds{\left[\sum \l_i e_i \right]}$ 
would be a nontrivial element of $\hker{E'}$,
contradicting $P(l)$ applied to the heap $E'$.  We conclude that $v = \be$.
The heap $E'(\be)$ inherits property P2 from $E'$ and is therefore
acyclic.  Since $E(\be) = \a \circ E'(\be)$ and $\hker{E(\be)} \ne 0$,
we apply statement $Q(l)$ to $E(\be)$ and conclude there is a minimal
vertex $\g$ of $E'(\be)$ with $\e(\g) = \e(\a)$.  Now $\a < \be < \g$ is
a convex chain in $E$ with $\e(\a) = \e(\g)$, which contradicts the fact
that $E$ has property P2 and completes the proof.
\qed\enddemo

\remark{Remark 2.4.3}
It is possible to find classes of heaps $H(P, \C)$ for which property P2
does not imply property P1.  For example, the heap in Remark 2.2.8 has
property P2 but not property P1.  In \S3.4, we shall look at some examples
of classes of heaps $H(P, \C)$ that do satisfy the hypotheses of Theorem
2.4.2.
\endremark

\proclaim{Theorem 2.4.4}
In a regular class of heaps, every acyclic heap is dismantlable, so 
property P1 and the property of being acyclic coincide in this case.
\endproclaim

\demo{Proof}
The second assertion is immediate from the first and Proposition 2.2.3.

Let $E = [E, \leq, \e]$ be an acyclic heap of the regular class $H(P, \C)$.  
If $E$ has property
P2, we are done.  If not, Remark 2.3.6 shows that there is a balanced convex
chain in $E$ of one of the types mentioned in lemmas 2.3.4 and 2.3.5.
Since $E$ is finite, there is a sequence $$
E = E_0, E_1, \ldots, E_l
$$ where $E_l$ has property P2 and for each $i$,
$E_{i+1}$ is obtained from $E_i$ by contraction of a balanced convex chain as
above.  Since $E$ is acyclic, and $E_l$ is acyclic
by Proposition 2.2.3, we see from lemmas 2.3.4 and 2.3.5 
that there will never be an opportunity to apply Lemma 2.3.4 in this
sequence.  We will be done if we can show that if $\bc$ is a balanced convex
chain $x < y < z$ for which $\e(x) \ne \e(y)$, then $E$ is dismantlable
if $E/\bc$ is dismantlable.

We proceed by induction on $n = |E/\bc|$.  If $n = 1$, the heap $E$ consists
solely of the chain $x < y < z$, which is dismantlable by inspection.  Suppose
the statement is true for $n = l$, and that $|E/\bc| = l + 1$.
Choose $\a$ such that $(E/\bc)(\a) \prec E/\bc$ as in Definition 2.2.1.
We will deal with the case $(E/\bc)(\a) \prec^- E/\bc$, the other case being
similar.

If $\a \ne x$ then the heap $E(\a)$ can be contracted
by Lemma 2.3.5 to the heap $(E/\bc)(\a)$; the latter heap has property
P1 by construction and $E(\a)$ is dismantlable by the inductive hypothesis.
Let $\be$ be minimal in $(E/\bc)(\a)$ but not in $E/\bc$, with $\e(\be) \ne
\e(\a)$.  Then $\be$ lies
in $E(\a)$ and is minimal in $E(\a)$ but not in $E$.  This shows that $E$ is
dismantlable.

The other possibility is that $\a = x$, in which case we have a sequence $$
E/\bc \prec E(x) \prec E
,$$ and the claim follows from Definition 2.2.1.
\qed\enddemo

\head 3. Quotients of heap monoid algebras \endhead

We show in \S3 that the dimensions of $\hker{E}$ of a heap
$E$ have a nice interpretation as the structure constants of a certain 
algebra associated to the heap $E$.  This fact is our main motivation 
in this paper.

\subhead 3.1 The heap monoid and some related structures \endsubhead

We now introduce the heap monoid associated to an arbitrary heap.  This
is naturally isomorphic to the commutation monoid (or ``free partially
abelian monoid'') appearing in the work
of Cartier and Foata \cite{{\bf 2}}.  The Mazurkiewicz traces 
\cite{{\bf 13}} used 
to study concurrency in computer science are another variant of the same idea,
and there is a large body of literature about them.

\definition{Definition 3.1.1}
A class of heaps $H(P, \C)$ has a natural monoid structure with composition
given by the superposition map $\circ$ of Definition 1.1.6.  We call this
monoid the {\it heap monoid}.
\enddefinition

\definition{Definition 3.1.2}
Let $A$ be a set and let $A^*$ be the free monoid generated by $A$.  Let $C$
be a symmetric and antireflexive relation on $A$.  The 
{\it commutation monoid} $\text{\rm Co} (A, C)$ is the quotient of the free
monoid $A^*$ by the congruence $\equiv_C$ generated by the commutation
relations: $$
ab \equiv_C ba \text{ for all } a, b \in A \text{ with } a \  C \ b
.$$
\enddefinition

To explain the relationship between heap monoids and commutation monoids, it 
is convenient to consider linear extensions of heaps.  These are also 
discussed in \cite{{\bf 17}, \S3} and \cite{{\bf 15}, \S1.2}.

\definition{Definition 3.1.3}
Let $(E, \leq)$ be a poset with $n$ elements.  A {\it natural labelling} of
$(E, \leq)$ is a bijection $\pi: E \ra [n] = \{1, 2, \ldots, n\}$ such that
$\a < \be$ implies that $\pi(\a) < \pi(\be)$.  If $[E, \leq, \e]$ is a
heap of $H(P, \C)$ and $p_i = \e(\pi^{-1} (i))$, we call the word $$
p_1 p_2 \cdots p_n
$$ of $P^*$ a {\it linear extension} of $E$.
\enddefinition

The following is a standard result about heaps, and a proof may be found in
\cite{{\bf 17}, Proposition 3.4}.

\proclaim{Proposition 3.1.4}
Let $E = [E, \leq, \e]$ be a heap of $H(P, \C)$, let $C$ be the 
complementary relation of $\C$, and let $\pi$ be a natural labelling
of $E$.  If we regard the words $P^*$ as representing elements of
$\text{\rm Co}(P, C)$, then the map sending $E$ to its linear extension
in $P^*$ under $\pi$ is independent of the choice of $\pi$,
and induces an isomorphism of monoids $H(P, \C) \ra \text{\rm Co}(P, C)$. \qed
\endproclaim

The following quotient of the monoid algebra will be of interest in our
applications.

\definition{Definition 3.1.5}
Maintain the above notation.  Let $\A$ be the ring of
Laurent polynomials $\zed[v, v^{-1}]$, let $\d := v + v^{-1}$, 
and let $\A \text{\rm Co}(P, C)$ be the monoid
algebra of $\text{\rm Co}(P, C)$ over $\A$.  We define the {\it generalized 
Temperley--Lieb algebra} $TL(P, \C)$ to be the $\A$-algebra 
obtained by quotienting $\A \text{\rm Co}(P, C)$ by the relations $$\eqalign{
s s &= \d s,\cr
s t s &= s \text{ if } s \ne t \text{ and } s \ \C \ t,\cr
}$$ where $s, t \in P$.
\enddefinition

\proclaim{Lemma 3.1.6}
The isomorphism of Proposition 3.1.4 induces an isomorphism between the 
algebra $TL(P, \C)$ and the quotient of $\A H(P, \C)$ by 
the relations $$\eqalign{
E &= \d E/\bc \text{ if } \bc \text{ is a balanced convex chain of 
length 2},\cr
E &= E/\bc \text{ if } \bc \text{ is a balanced convex chain } x < y < z
\text{ with } \e(x) \ne \e(y).\cr
}$$
\endproclaim

\demo{Proof}
If $p_1 p_2 \cdots p_n$ is a word in $P^*$ corresponding to a heap 
$E = [E, \leq, \e]$ in $H(P, \C)$ with natural labelling $\pi$, then equation
(7) of \cite{{\bf 17}, \S3} shows that $\pi^{-1}(i) < \pi^{-1}(j)$ 
if and only if there is a sequence $$
1 \leq i = i_1 < \cdots < i_t = j \leq n
$$ such that $p_{i_m} \ \C \ p_{i_{m+1}}$ for $1 \leq m < t$.
It follows from this observation that subwords of 
$p_1 p_2 \cdots p_n$ of the form $ss$
(respectively, $sts$) as in Definition 3.1.5 correspond to chains in $E$ of 
length $2$ (respectively, $3$) as described in the statement.

The converse implication follows from 
the standard fact that for any convex chain in a poset, there
exists a linear extension of the poset in which the members of the chain
appear consecutively.
\qed\enddemo

\subhead 3.2 Structure constants for $TL(P, \C)$ \endsubhead

In \S3.2, we exhibit a free $\A$-basis for $TL(P, \C)$ using Bergman's
diamond lemma \cite{{\bf 1}}, and show that in favourable circumstances the
structure constants of the algebra with respect to this basis are closely
related to dimensions of $\ker \pd$ for certain heaps.  This allows us
to explore our main application.  We remark that the idea for finding this
basis essentially comes from Graham's thesis \cite{{\bf 7}, Theorem 6.2}.

In order to use Bergman's diamond lemma we need to recall some terminology
from \cite{{\bf 1}}.  Let $R$ be a commutative ring and let $X$ be a nonempty
set.  Let $\leq_X$ be a semigroup partial order on $X^*$: that is, if
$\l, \mu, \nu$ are (possibly empty) words in $X^*$ and $\mu \leq \nu$, then
we have $\l \mu \leq \l \nu$ and $\mu \l \leq \nu \l$.  We say $\leq_X$
satisfies the descending chain condition if any sequence $\l_1 >_X \l_2 >_X 
\cdots$ terminates.  A reduction system $S$ for $R X^*$ is a set of rules 
of the form $s : \mu_s \mapsto a_s$, where $\mu_s$ lies in $X^*$ and $a_s$ 
lies in $R X^*$.  The $R$-module 
maps $R X^* \ra R X^*$ used to apply rules are known as reductions;
these may consist of several rules performed sequentially.
The two-sided ideal ${\Cal I}(S)$ of $R X^*$
is that generated by all elements $\mu_s
- a_s$ for all rules $s \in S$.  We say $s$ is {\it compatible} 
with $\leq_X$ if 
$a_s$ can be written as a linear combination of monomials strictly 
less than $\mu_s$ in $\leq_X$, and we say $S$ is compatible with $\leq_X$ if
each of its rules is.

An {\it overlap ambiguity} occurs when there are two rules $s_1$ and $s_2$
such that there exist monomials $\nu_2$ and $\nu_1$ with $\mu_{s_1} \nu_2
= \nu_1 \mu_{s_2}$; it is said to be {\it resolvable} if there are reductions
$t_1$ and $t_2$ with $t_1(a_{s_1} \nu_2) = t_2(\nu_1 a_{s_2})$.  An 
{\it inclusion ambiguity} occurs when there are two rules $s_1$ and $s_2$
such that there exist monomials $\l$ and $\nu$ with $\mu_{s_2} = \l \mu_{s_1}
\nu$; it is said to be {\it resolvable} if there are reductions
$t_1$ and $t_2$ with $t_1(\l a_{s_1} \nu) = t_2(a_{s_2})$.

A reduction $t$ is said to act trivially on $a \in R X^*$ if $t(a) = a$, and
if all reductions act trivially on $a$, we say $a$ is irreducible.  The set of
irreducible elements arising from $S$ is denoted $\text{\rm Irr} (S)$.  A
normal form of $a \in R X^*$ is an element $b \in \text{\rm Irr} (S)$ to which
$a$ can be reduced; it is not immediate that normal forms always exist or
that they are unique.

The following theorem is part of Bergman's diamond lemma, which
is proved in \cite{{\bf 1}, Theorem 1.2}.

\proclaim{Theorem 3.2.1 (Bergman)}
Let $R$ be an associative, commutative ring with $1$.  Let $X$ be a nonempty
set, let $\leq_X$ be a semigroup partial order on $X^*$ and let $S$ be a
reduction system for $R X^*$.  If $S$ is compatible with $\leq_X$ and $\leq_X$
satisfies the descending chain condition then the following are equivalent:
\item{\rm (i)}{All ambiguities in $S$ are resolvable.}
\item{\rm (ii)}{Every element $a \in R X^*$ has a unique normal form which 
equals $t(a)$ for some reduction $t$.}
\item{\rm (iii)}{$R X^* = \text{\rm Irr} (S) \oplus {\Cal I}(S)$.}
\qed\endproclaim

\proclaim{Proposition 3.2.2}
The quotient of $\A H(P, \C)$ described in Lemma 3.1.6 has as a free
$\A$-basis the images of those heaps in $H(P, \C)$ with property P2.
\endproclaim

\demo{Proof}
We use Theorem 3.2.1.  Let $R = \A = \zed[v, v^{-1}]$, let $X = P$ and let 
$\leq_X$ be the partial order such that $b_1 b_2 \cdots b_l \leq_X
c_1 c_2 \cdots c_m$ if and only if $b_1 b_2 \cdots b_l$ is a subsequence of
$c_1 c_2 \cdots c_m$; this is a semigroup partial
order.  There are two kinds of reduction rules.

For the first kind of reduction rule, we take $$
\mu_1 = p p_{i_1} \cdots p_{i_l} p
$$ where the $p_{i_n}$ (if there are any) are distinct from $p$ and
commute with $p$.  This is precisely
the condition for the occurrences of $p$ to 
correspond by Proposition 3.1.4 to heap elements $\a$
and $\be$ for which $\a < \be$ is a balanced convex chain.
The element $a_1$ in this case is given by $$
\d p p_{i_1} \cdots p_{i_l} = 
\d p_{i_1} \cdots p_{i_l} p.$$

For the second kind of reduction rule, we take $$
\mu_2 = p p_{i_1} \cdots p_{i_l} p' p_{j_1} \cdots p_{j_m} p
$$ where all the $p_{i_n}$ and $p_{j_n}$ (if there are any) are distinct
from $p$ and commute with $p$, and $p'$ does not commute with $p$.  This is
precisely the condition for the letters $p$, $p'$, $p$ to 
correspond respectively to heap elements $\a$, $\be$, $\g$
for which $\a < \be < \g$ is a balanced convex chain with
$\e(\a) \ne \e(\be)$.  The element $a_2$ in this case is given by $$
a_2 = p p_{i_1} \cdots p_{i_l} p_{j_1} \cdots p_{j_m}
= p_{i_1} \cdots p_{i_l} p_{j_1} \cdots p_{j_m} p
.$$

In each case, $a_i$ is a multiple of a strictly shorter monomial than $\mu_i$.
Since every heap is a finite set, $\leq_X$ has the descending chain
condition.

We now show that all ambiguities are resolvable.  Most of the possible
inclusion ambiguities are easily seen to be resolvable.  The only difficult
case arises from words of the form $$
p \cdots (p'' \cdots p' \cdots p'') \cdots p
,$$ where $p$ and $p''$ are distinct and commute with each other but neither 
commutes with $p'$, and where the both the whole word and the 
parenthetic expression shown are of the form $\mu_2$ as above.  This type of
ambiguity resolves to $$
\d p \cdots (p'' \cdots \widehat{p'} \cdots \widehat{p''}) \cdots \widehat{p}
$$ (where the hats denote omission) 
whether the outermost or the innermost reduction rule is applied first.

The overlap ambiguities correspond to chains in the heap
of the following kinds:

\item{1.}{Convex chains $b < c < d$ with $\e(b) = \e(c) = \e(d)$.}
\item{2.}{Convex chains $b < c < d < e$ with either 
$\e(b) = \e(c) = \e(e) \ne \e(d)$
or
$\e(b) = \e(d) = \e(e) \ne \e(c)$ .
\item{3.}{Chains $b < c < d < e < f$ with $b < c < d$ convex and
$d < e < f$ convex and $\e(c) \ne \e(b) = \e(d) = \e(f) \ne \e(e)$ (but
not necessarily with $\e(c) = \e(e)$).}
}

These ambiguities are easily seen to be resolvable.  The definition of
property P2 shows that a word in $P^*$ will be irreducible in this reduction
system if and only if it has property P2, so by Theorem 3.2.1 the heaps
with property P2 represent a basis for the quotient algebra.
\qed\enddemo

The following theorem gives a nice interpretation of the dimensions of
$\ker \pd$ of heaps in terms of the algebra $TL(P, \C)$ and its structure
constants.

\proclaim{Theorem 3.2.3}
\item{\rm (i)}{Let $D$ be an arbitrary heap of $H(P, \C)$.  Then there exists
a unique heap $G \in H(P, \C)$ with property P2 such that $D = \d^m G$ in the
quotient algebra $TL(P, \C)$.  We have $$\dim \hker{D} = m + \dim \hker{G}.$$
In particular, if $H(P, \C)$ is a regular class of heaps then 
$m = \dim \hker{D}$.}
\item{\rm (ii)}
{Suppose $H(P, \C)$ is a regular class of heaps and let $E$ and $F$ be 
heaps with property $P2$ regarded as basis elements of $TL(P, \C)$ via 
Proposition 3.2.2.  Then there is a basis element $G$ of $TL(P, \C)$ for
which $E \circ F = \d^{\dim \hker{E \circ F}} G$.}
\endproclaim

\demo{Proof}
Part (ii) follows from part (i), so we prove the former.

The existence of $G$ comes from Theorem 3.2.1 and Proposition 3.2.2: the
normal form of $D$ will be a multiple of only one basis element
because the reduction rules of Proposition 3.2.2 take monomials to multiples
of monomials.  The equation $$
\dim \hker{D} = m + \dim \hker{G}
$$ comes from comparing Lemma 3.1.6 with lemmas 2.3.4 and 2.3.5.
We apply the latter two lemmas repeatedly to the heap $D$ until
no further reductions are possible and we are left with the heap $G$.  Along
the way, the dimension of $\ker \pd$ will decrease by $1$
precisely when an extra factor of $\d$ appears in Lemma 3.1.6.  For the
last claim, observe that if $H(P, \C)$ is regular then $G$ has property P1
because it has property P2, and the conclusion follows from Proposition 2.2.3.
\qed\enddemo

\subhead 3.3 Heaps arising from Coxeter groups \endsubhead

In \S3.4, we will apply the theory developed here to results involving Coxeter
groups, and for this we need to explain the connection between Coxeter
groups and heaps.  This is a major theme of the paper \cite{{\bf 15}}.

We only consider the case of simply laced Coxeter groups, as this is all we
need for our purposes.  This simplifies the definitions somewhat.

\definition{Definition 3.3.1}
Let $\Gamma$ be a graph.  The (simply laced) Coxeter group $W(\Gamma)$ 
associated to $\Gamma$
is given by generators $\{s_i : i \in S\}$ where $S$ is the set of vertices
of $\Gamma$, and defining relations $$\eqalign{
s_i^2 &= 1,\cr
s_i s_j s_i &= s_j s_i s_j \quad 
\text{ if } i, j \text{ are adjacent in } \Gamma,\cr
s_i s_j &= s_j s_i \quad \text{ otherwise.}\cr
}$$
\enddefinition

We now define the fully commutative elements of a Coxeter group.  
These were studied in Fan's thesis \cite{{\bf 4}} in the simply laced case under
the name ``commutative elements'', and a general definition was given by
Stembridge \cite{{\bf 15}, \S1.1}.

\definition{Definition 3.3.2}
Let $W = W(\Gamma)$ be a Coxeter group as above.  Clearly every element 
$w$ can be written as $$
w = s_{i_1} s_{i_2} \cdots s_{i_r}
$$ for some generators $s_i$.  If $r$ is minimal for a given $w$, we define the
{\it length} of $w$ to be $r$, and we call the associated product of generators
{\it reduced}.  If any two reduced expressions for $w$ may be transformed
into each other by repeated applications of relations of the form  
$s_i s_j = s_j s_i$ as in 
Definition 3.3.1, we say $w$ is {\it fully commutative}.  The set of
fully commutative elements of $W$ is denoted by $W_c$.
\enddefinition

Elements of Coxeter groups give rise to heaps as follows; 
see also \cite{{\bf 15},
\S1.2}.  The next procedure is inverse to that given by Proposition 3.1.4.

\definition{Definition 3.3.3}
Let $\Gamma$ be a graph, let $S$ be the set of vertices of $\Gamma$ and let
$\C$ be the relation on $S$ defined by $s_i \ \C \ s_j$ if and only if either
$s_i = s_j$ or $s_i$ and $s_j$ are adjacent vertices.  Let $C$ be the 
complementary relation of $\C$.  Let $w = s_{i_1} \cdots s_{i_l} \in S^*$
be an arbitrary word in the generators $S$
of a simply laced Coxeter group $W = W(\Gamma)$.  The word
$w$ gives a labelled heap $(E, \leq_E, \e)$ where $E = \{1, 2, \ldots, l\}$,
$\e(j) = s_{i_j}$, and the relation $\leq_\C$ of condition 2 of 
Definition 1.1.1 can be defined by $$
\a \leq_\C \be \Leftrightarrow \a \leq \be \text{ and } \e(\a) \ \C \ \e(\be)
,$$  where $\leq$ is the usual ordering on integers.
The partial order $\leq_E$ is the transitive extension of $\leq_\C$,
and the heap of $H(S, \C)$ 
corresponding to the given labelled heap is by definition the heap
of the word $w$.  Moreover, the heap of $w$ may be recovered from the element
of $\text{\rm Co} (S, C)$ represented by $w$.
\enddefinition

The following result is a special case of Stembridge's 
\cite{{\bf 15}, Proposition 2.3}.

\proclaim{Proposition 3.3.4}
Let $\Gamma$ be a graph with set of vertices $S$.  A word in $S^*$ is a
reduced expression for 
an element of $W_c \subset W = W(\Gamma)$ if and only if the associated
heap (as in Definition 3.3.3) has property P2.
\qed\endproclaim

We can use the graph $\Gamma$ to define an analogue of the algebra
$TL(P, \C)$ of \S3.1.

\definition{Definition 3.3.5}
Let $\Gamma$ be a graph with set of vertices $S$ 
and let $\A = \zed[v, v^{-1}]$ with $\d = v + v^{-1}$.  We define the
$\A$-algebra $TL(\Gamma)$ by generators $\{e_i : i \in S\}$ and 
relations $$\eqalign{
e_i^2 &= \d e_i,\cr
e_i e_j e_i &= e_i \quad 
\text{ if } i, j \text{ are adjacent in } \Gamma,\cr
e_i e_j &= e_j e_i \quad \text{ otherwise.}\cr
}$$
\enddefinition

\remark{Remark 3.3.6}
Definition 3.3.5 is compatible with Definition 3.1.5: if $H(P, \C)$ is a class
of heaps with concurrency graph $\Gamma$, then $TL(\Gamma)$ is canonically
isomorphic to $TL(P, \C)$.
\endremark

The following result is well known.

\proclaim{Proposition 3.3.7}
Let $\Gamma$ be a graph with vertex set $S$ and let $w$ be a fully commutative
element of the corresponding Coxeter group with reduced expression $
s_{i_1} s_{i_2} \cdots s_{i_r}
$.  Then the element $e_w := e_{i_1} e_{i_2} \cdots e_{i_r}$ is independent of
the choice of reduced expression and the set $\{e_w : w \in W_c\}$ is a
free $\A$-basis of $TL(\Gamma)$ (called the {\it monomial basis}).
\endproclaim

\demo{Proof}
The first assertion is an easy consequence of the definition of full
commutativity.  The second assertion follows from Proposition 3.2.2, Lemma
3.1.6 and Proposition 3.3.4.
\qed\enddemo

\subhead 3.4 An application to algebra \endsubhead

We will be particularly interested in heaps whose concurrency graph is of
type $E_n$, namely the graph shown in Figure 2.

\topcaption{Figure 2} Coxeter graph of type $E_n$
\endcaption
\centerline{
\hbox to 1.500in{
\vbox to 0.500in{\vfill
        \includegraphics{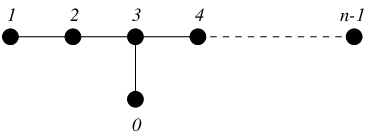}
}
\hfill}
}

The numbering above is chosen to agree with that in Fan's paper \cite{{\bf 5},
\S6.3}.  The elements of $P$ will be the integers $0, 1, 2, \ldots, n-1$.
We emphasize that $n$ is arbitrary; in particular, it can be greater than $8$.
Graphs of type $A_{n-1}$ (respectively, $D_{n-1}$) may be obtained from the 
graph in Figure 2 by omitting the vertices numbered $0$ (respectively, 1).
Graphs that are isomorphic to graphs of type $A_n$, $D_n$ or $E_n$ will be 
called {\it graphs of type $ADE$}.  We will also be concerned with
graphs of type $\widetilde{A}_{n-1}$ for $n \geq 3$; these are $n$-gons
with vertices numbered consecutively $1$ up to $n$.

Interesting examples of our main results arise as follows.

\vfill\eject
\proclaim{Theorem 3.4.1 (Fan)}
Let $H(P, \C)$ be a class of heaps whose concurrency graph is of type $ADE$
or of type $\widetilde{A}_{n-1}$ with $n$ odd.
Then $H(P, \C)$ is a regular class of heaps.
\endproclaim

\demo{Proof}
The $ADE$ case is a restatement of \cite{{\bf 5}, Lemma 4.3.1}, and the case of
type $\widetilde{A}_{n-1}$ is a restatement of \cite{{\bf 6}, 
Proposition 3.1.2}.
\qed\enddemo

The following subtle property of the algebras $TL(\Gamma)$ (see Definition
3.3.5) is originally due in the $ADE$ case 
to Graham \cite{{\bf 7}, Lemma 9.10, Lemma 9.13}, who proved it
using a complicated combinatorial argument.  It is of key importance in
certain applications as we explain in \S4.

\proclaim{Proposition 3.4.2}
Let $\Gamma$ be the concurrency graph of a regular class of heaps, $H(P, \C)$.
Consider an arbitrary word $$
e = e_{i_1} e_{i_2} \cdots e_{i_r}
$$ in the generators for $TL(\Gamma)$, and define $w \in W_c$ such that  
$e = \d^m e_w$ for $w \in W_c$.
\item{\rm (i)}{If we omit one generator from the word $e$ to form $$
e' = e_{i_1} \cdots \widehat{e_{i_l}} \cdots e_{i_r}
$$ (where the hat denotes omission) and write $e' = \d^{m'} e_{w'}$
for some $w' \in W_c$, then $m' \in \{m-1, m, m+1\}$.}
\item{\rm (ii)}{If $e$ is of the form $e_x$ for some $x \in W_c$ then
$m' = 0$.}
\endproclaim

\demo{Proof}
Part (i) is a consequence of the Deletion Lemma (Theorem 2.1.1) and Theorem
3.2.3 (i), first with $D = e, G = e_w$, and then with $D = e', G = e_{w'}$.

For part (ii), Theorem 2.4.2 shows that the heap of $x \in W_c$ is
strongly acyclic, and the conclusion follows from Theorem 3.2.3 (i).
\qed\enddemo

\head 4. Applications and questions \endhead

We conclude with a survey of how the results of this paper are related to
results in the literature, and some possible directions for future research.

\subhead 4.1 Canonical bases for Hecke algebra quotients \endsubhead

In \cite{{\bf 9}}, the author and J. Losonczy introduced canonical bases
(IC bases) for the generalized Temperley--Lieb algebras of \cite{{\bf 7}} and
showed that in the case of the algebras $TL(\Gamma)$ where $\Gamma$ is of
type $ADE$, the basis of Proposition 3.2.2 is the canonical basis.  This 
relies heavily on Proposition 3.4.2.  An argument similar to our 
proof of Proposition 3.4.2
shows that any strongly acyclic heap corresponds to a canonical basis element
for $TL(\Gamma)$ that is given by a monomial in the generators $e_i$.  (Not
all canonical basis elements are of this form.)

Although Theorem 2.4.2 (ii) has so far only appeared as a by-product, it
has some nice applications of its own, one of which is to allow a recurrence
formula for canonical basis elements for $TL(\Gamma)$ similar to that given
by Kazhdan and Lusztig \cite{{\bf 11}, \S2.2} in the case of Hecke algebras.
It also gives recurrence formulae for analogues of inverse Kazhdan--Lusztig
polynomials (denoted by $\widetilde{Q}_{x, w}$ in \cite{{\bf 9}, Lemma 3.5}).

Our results also have applications to generalized Temperley--Lieb
algebras arising from non-simply-laced Coxeter graphs.  These algebras
have bases indexed by heaps satisfying conditions similar to, but weaker 
than property P2; see \cite{{\bf 15}, Proposition 2.3} for the exact condition.
The canonical bases for these algebras are described in \cite{{\bf 8}} for 
Coxeter types $B$ and $H$.  In these two special cases, the heaps
indexing the basis are acyclic, which is essentially a neater restatement of 
\cite{{\bf 8}, Lemma 3.1.1}; the Deletion Lemma (Theorem 2.1.1) then gives
\cite{{\bf 8}, Lemma 3.1.3}.  Another result of \cite{{\bf 8}} that is 
more easily
phrased in terms of the map $\pd$ is \cite{{\bf 8}, Proposition 3.1.9}, which 
classifies the possible image vertices in a basis heap.

In some cases, such as Coxeter type $B$, it is possible to find for 
every element $w$ in the Coxeter group a reduced expression whose heap 
(see Definition 3.3.3) is acyclic.  The results of \cite{{\bf 10}, \S2} are
essentially consequences of this observation.  The analogous claim
in type $D$ is false, which makes Losonczy's argument in
\cite{{\bf 12}} much more difficult.

We remark that if the basis heaps mentioned above are acyclic, it is possible
to state a (more complicated) version of Theorem 3.2.3 for the 
non-simply-laced case.

We hope to give details of these applications separately.

\subhead 4.2 Computing $\dim \hker{E}$ \endsubhead

Certain of the algebras $TL(\Gamma)$ may be understood by a calculus of
diagrams rather than using generators and relations.  In these cases, a
word in the generators $e_i$ of $TL(\Gamma)$ may be represented by a diagram,
and Proposition 3.1.4 then shows how to represent a heap as a diagram.
We outline here how these diagrams may be used to calculate $\dim
\hker{E}$ essentially by inspection.
The proofs are too long to present here, particularly in
the case of type $E_n$ below, but we hope to give details separately.

\subsubhead Type $A_n$ \endsubsubhead

Suppose $H(P, \C)$ is a class of heaps with concurrency graph $\Gamma$ of
type $A_n$, in other words, $\Gamma$ is a line.  In this case, $TL(\Gamma)$
is the Temperley--Lieb algebra of \cite{{\bf 16}} which has a well known diagram
calculus whose origins can be traced back to \cite{{\bf 14}}.  If $E$ is a heap
of $H(P, \C)$ then $\dim \hker{E}$ is the number of loops occurring in the
corresponding diagram.

\subsubhead Type $E_n$ \endsubsubhead

Suppose $H(P, \C)$ is a class of heaps with concurrency graph $\Gamma$ of
type $E_n$, so that it is isomorphic to a graph such as
that shown in Figure 2.  Let $E$ be a heap in a class
$H(P, \C)$ that has concurrency graph $\Gamma$.  A diagram calculus for 
$TL(\Gamma)$ was described by tom Dieck in \cite{{\bf 3}}, and in fact 
this gives a
faithful representation of the algebra, although this is not proved in 
\cite{{\bf 3}}.  In this case the diagrams are non-intersecting curves drawn 
inside a rectangle whose endpoints lie on the boundary, and there are certain
discs (``pillars'') lying in the connected components of the complement of
these curves.  If $E$ is a heap of $H(P, \C)$ then $$
\dim \hker{E} = a + \sum_{i \in I} \max(0, \ b_i - 1)
,$$ where $a$ is the number of loops containing no pillars, $I$ is the set of
connected components in the diagram (including the insides of loops) 
and $b_i$ is the number of pillars in component $i$.

\subsubhead Type $\widetilde{A}_{n-1}$ \endsubsubhead

Suppose $H(P, \C)$ is a class of heaps with concurrency graph $\Gamma$ of
type $\widetilde{A}_{n-1}$ ($n \geq 3$); in other words, let $\Gamma$ be an 
$n$-gon.  A diagram calculus for $TL(\Gamma)$ is given in \cite{{\bf 6}, \S4}.
The diagrams consist of certain curves inscribed on the surface of a cylinder
whose endpoints (if any) lie on the boundary.  
A diagram contains a number $a$ of 
loops contractible on the cylinder, and a number $b$ of loops that are not
contractible on the cylinder.  (These numbers may be zero.)  If $E$ is a 
heap of $H(P, \C)$ then $$
\dim \hker{E} = a + c (\max(0, \ b - 1))
,$$ where $c = 1$ if $4 | n$ or $k$ has characteristic $2$, and $c = 0$ 
otherwise. 

\subhead 4.3 Concluding questions \endsubhead

In light of the results of this paper, it would be interesting to have an
answer to the following graph theoretic problem.

\proclaim{Problem 4.3.1}
Find necessary and sufficient conditions on a graph $\Gamma$ for it to
be the concurrency graph of a regular class of heaps.
\endproclaim

In a future paper, we will present a complete solution to Problem 4.3.1
under the assumption that the set of pieces, $P$, is finite.  That paper
will also show how, with considerably more work, Theorem 2.4.2 may 
be sharpened.

Another intriguing direction for future research is suggested by the following

\proclaim{Question 4.3.2}
Do the results of this paper have applications to concurrency in computer
science?
\endproclaim

\head Acknowledgements \endhead

I wish to thank J. Losonczy for helpful comments,
and Colorado State University for its hospitality during the
preparation of this paper.

\leftheadtext{}
\rightheadtext{}

\end
\Refs\refstyle{A}\widestnumber\key{{\bf 17}}
\leftheadtext{References}
\rightheadtext{References}

\ref\key{{\bf 1}}
\by G.M. Bergman
\paper The diamond lemma for ring theory
\jour Adv. Math. \vol 29 \yr 1978 \pages 178--218
\endref

\ref\key{{\bf 2}}
\by P. Cartier and D. Foata
\paper Probl\`emes combinatoires de commutation et r\'earrangements
\jour Lecture Notes in Mathematics
\vol 85
\yr 1969
\publ Springer-Verlag
\publaddr New York/Berlin
\endref

\ref\key{{\bf 3}}
\by T. tom Dieck
\paper Bridges with pillars: a graphical calculus of knot algebra
\jour Topology Appl.
\vol 78 \yr 1997 \pages 21--38
\endref

\ref\key{{\bf 4}}
\by C.K. Fan
\book A Hecke algebra quotient and properties of commutative elements
of a Weyl group
\publ Ph.D. thesis
\publaddr M.I.T.
\yr 1995
\endref

\ref\key{{\bf 5}}
\by C.K. Fan
\paper Structure of a Hecke algebra quotient
\jour J. Amer. Math. Soc.
\vol 10  \yr 1997 \pages 139--167
\endref

\ref\key{{\bf 6}}
\by C.K. Fan and R.M. Green
\paper On the affine Temperley--Lieb algebras
\jour Jour. L.M.S.
\vol 60 \yr 1999 \pages 366--380
\endref

\ref\key{{\bf 7}}
\by J.J. Graham
\book Modular representations of Hecke algebras and related algebras
\publ Ph.D. thesis
\publaddr University of Sydney
\yr 1995
\endref

\ref\key{{\bf 8}}
\by R.M. Green
\paper Decorated tangles and canonical bases
\jour J. Algebra
\vol 246 \yr 2001 \pages 594--628
\endref

\ref\key{{\bf 9}}
\by R.M. Green and J. Losonczy
\paper Canonical bases for Hecke algebra quotients
\jour Math. Res. Lett.
\vol 6 \yr 1999 \pages 213--222
\endref

\ref\key{{\bf 10}}
\by R.M. Green and J. Losonczy
\paper A projection property for Kazhdan--Lusztig bases
\jour Int. Math. Res. Not.
\vol 1 \yr 2000 \pages 23--34
\endref

\ref\key{{\bf 11}}
\by D. Kazhdan and G. Lusztig
\paper Representations of Coxeter groups and Hecke algebras
\jour Invent. Math.
\vol 53 \yr 1979 \pages 165--184
\endref

\ref\key{{\bf 12}} \by J. Losonczy \paper The Kazhdan--Lusztig
basis and the Temperley--Lieb quotient in type D \jour J. Algebra
\vol 233 \yr 2000 \pages 1--15
\endref

\ref\key{{\bf 13}}
\by A. Mazurkiewicz
\paper Trace theory
\inbook Petri nets, applications and relationship to other models of
concurrency
\publ Lecture Notes in Computer Science 255, Springer
\publaddr Berlin--Heidelberg--New York
\yr 1987
\pages 279--324
\endref

\ref\key{{\bf 14}}
\by R. Penrose
\paper Angular momentum: an approach to combinatorial space-time
\inbook Quantum Theory and Beyond (E. Bastin, Ed.)
\publ Cambridge University Press
\publaddr Cambridge
\yr 1971
\pages 151--180
\endref

\ref\key{{\bf 15}}
\by J.R. Stembridge
\paper On the fully commutative elements of Coxeter groups
\jour J. Algebraic Combin.
\vol 5
\yr 1996
\pages 353--385
\endref

\ref\key{{\bf 16}}
\by H.N.V. Temperley and E.H. Lieb
\paper Relations between percolation
and colouring problems and other graph theoretical problems associated
with regular planar lattices: some exact results for the percolation
problem
\jour Proc. Roy. Soc. London Ser. A
\vol 322 \yr 1971 \pages 251--280
\endref

\ref\key{{\bf 17}}
\by G.X. Viennot
\paper Heaps of pieces, I: basic definitions and combinatorial lemmas
\inbook Combinatoire \'E\-nu\-m\'e\-ra\-tive
\publ Springer-Verlag
\publaddr Berlin
\yr 1986 \pages 321--350 \bookinfo ed. G. Labelle and P. Leroux
\endref


\endRefs
